\documentclass{CSML}

\pdfoutput=1

\usepackage{lastpage}

\def\dOi{13(4:24)2017}
\lmcsheading%
{\dOi}
{1--\pageref{LastPage}}
{}
{}
{Jul.~\phantom06, 2014}
{Dec.~\phantom07, 2017}
{}

\usepackage{amsmath,amssymb,amsfonts,latexsym,stmaryrd,mathrsfs}
\usepackage{cmll,bbm} 
\usepackage[hidelinks]{hyperref}

\newcommand{\pair}[2]{\langle #1 , #2 \rangle}
\newcommand{\arr}{\rightarrow}
\newcommand{\Pow}{\mathcal{P}}
\newcommand{\Set}{{\mathbf{Set}}}

\newcommand{\Pfin}{\Pow_{\mathsf{fin}}}
\newcommand{\fin}{{\mathsf{fin}}}

\newcommand{\N}{\mathbb{N}}

\def\lbr{\mathopen{{[\kern-0.14em[}}}   
\def\rbr{\mathclose{{]\kern-0.14em]}}}  
\newcommand{\sem}[1]{\lbr#1\rbr}

\newcommand{\inc}{\hookrightarrow}

\newcommand{\strarr}{\mbox{$\circ \kern-0.4em \arr$}}
\newcommand{\eps}{\,\varepsilon\,}
\newcommand{\neps}{\; \slash \kern-0.7em \eps \,}

\newcommand{\lolli}{{\multimap}}

\newcommand{\RT}{{\mathbf{RT}}}

\newcommand{\BR}{{\mathsf{BR}}}

\newcommand{\C}{{\mathcal{C}}}

\newcommand{\E}{{\mathcal{E}}}
\newcommand{\F}{{\mathcal{F}}}

\newcommand{\id}{{\mathrm{id}}}

\newcommand{\forget}[1]{}

\newcommand{\mono}{\rightarrowtail}

\newcommand{\Sh}{\mathsf{Sh}}
\newcommand{\Se}{\mathbf{Set}}

\newcommand{\Mo}{{\mathcal{M}}}

\newcommand{\Eff}{\mathcal{E}\kern-0.14em\mathit{ff}}

\newcommand{\Asm}{{\mathbf{Asm}}}

\newcommand{\Kcal}{{\mathcal{K}}}

\newcommand{\schloop}{{< \kern-0.40em  <}}

\newcommand{\eq}{{\mathit{eq}}}
\newcommand{\Cc}{{\mathsf{cc}}}
\newcommand{\kk}{{\mathsf{k}}}
\newcommand{\ic}{{\mathsf{i}}}
\newcommand{\PL}{\mathsf{PL}}
\newcommand{\trace}{{\mathsf{tr}}}

\begin{document}

\title[Classical Realizability Models arising from Models of Untyped Lambda Calculus]{A Classical Realizability Model arising from \\
       a Stable Model of Untyped Lambda Calculus}

\author[T.~Streicher]{Thomas STREICHER}	
\address{Fachbereich 4 Mathematik TU Darmstadt\\
         Schlo{\ss}gartenstr.~7, D-64289, Germany}	
\email{streicher@mathematik.tu-darmstadt.de}  





\begin{abstract}
  \noindent 
In \cite{reusstreicher} it has been shown that $\lambda$-calculus with control
can be interpreted in any domain $D$ which is isomorphic to the domain of 
functions from $D^\omega$ to the $2$-element (Sierpi\'nski) lattice $\Sigma$.
By a theorem of A.~Pitts there exists a unique subset $P$ of $D$ such that
$f \in P$ iff $f(\vec{d}) = \bot$ for all $\vec{d} \in P^\omega$. The domain
$D$ gives rise to a \emph{realizability structure} in the sense of 
\cite{kriv2011} where the set of proof-like terms is given by $P$.

When working in Scott domains the ensuing realizability model coincides 
with the ground model $\Se$ but when taking $D$ within coherence spaces
we obtain a classical realizability model of set theory different from 
any forcing model. We will show that this model validates countable and 
dependent choice since an appropriate form of bar recursion is available 
in stable domains.
\end{abstract}

\maketitle
\begin{center} \emph{Dedicated to Pierre-Louis Curien at the occasion of his 60th Birthday}
\end{center}

\section*{Introduction}

In the first decade of this millenium J.-L.~Krivine has developed his theory
of \emph{classical realizability}, see e.g.\ \cite{kriv2009,kriv2011}, 
for higher order logic and set theory. Whereas intuitionistic realizability
is based on the notion of a \emph{partial combinatory algebra} (pca) classical
realizability is based on a notion of \emph{realizability algebra} as defined 
in \cite{kriv2011}. Both notions are incomparable since not every pca can be
extended to a realizability algebra and there are realizability algebras which
do not contain a pca as a substructure. Accordingly, not all classical 
realizability models appear as booleanizations of intuitionistic realizability
models as studied in \cite{jaapbook}.

In the current paper, however, we concentrate on a particular classical 
realizability model which appears as a boolean subtopos of a \emph{relative
realizability} topos (see \cite{jaapbook}). The starting point for this model
is the observation from \cite{reusstreicher} that the recursive domain
$D \cong \Sigma^{D^\omega}$ gives rise to a model for $\lambda$-calculus 
with control.
(Here $\Sigma = \{\bot,\top\}$ is the 2-element Sierpi\'nski lattice and 
$D^\omega$ is the countable product of $D$.) Since $D$ is a model of untyped
$\lambda$-calculus it is in particular a pca. By a theorem of A.~Pitts 
\cite{pitts96} there exists a unique subset $P$ of $D$ such that $t \in P$ 
iff $t(\vec{s}) = \bot$ for all $\vec{s} \in P^\omega$. Obviously, this
subset $P$ forms a sub-pca of $D$ thus giving rise to the relative 
realizability topos $\E = \RT(D,P)$ as described in \cite{jaapbook}. 
Notice that $\top_D \in D \setminus P$ and thus $U = \{\top_D\}$ gives rise
to a nontrivial truth value in $\E$ different from both $\top_\E$ and $\bot_\E$.
This $U$ (like any subterminal object of $\E$) induces a closure operator 
(\emph{aka} Lawvere-Tierney topology) $j_U(p) = (p \to U) \to U$ on $\E$. As
is well known the subtopos $\E_U$ of $j_U$-sheaves of $\E$ is boolean.

We will show that $\E_U$ is equivalent to the classical realizability topos
$\Kcal$ induced by the realizability structure whose set $\Lambda$ of terms 
is $D$, whose set of stacks $\Pi$ is $D^\omega$ and whose set $\PL$ of 
proof-like terms is $P$. We will show that $\Kcal$ is equivalent to $\Set$ 
when $D$ is the bifree solution of the domain equation 
$D \cong \Sigma^{D^\omega}$ in Scott domains. However, when considering 
the solution of $D \cong \Sigma^{D^\omega}$ in the category $\mathbf{Coh}$ of 
coherence spaces and Scott continuous and stable maps then the ensuing boolean
topos $\Kcal$ is not a Grothendieck topos and thus \emph{a fortiori} not a 
forcing model\footnote{i.e.\ a category of sheaves over a complete boolean 
algebra or, equivalently, a Grothendieck topos where every epimorphism splits,
see e.g.\ \cite{Elephant}}. 
We will show that $\Kcal$ validates all true sentences of first order 
arithmetic and the principles of countable and dependent choice.

\section{Realizability structures induced 
             by $D \cong \Sigma^{D^\omega}$}\label{rsd}

Quite generally we might consider objects $D \cong \Sigma^{D^\omega}$ in well
pointed cartesian closed categories $\C$ with countable products and an 
object $\Sigma$ having precisely two global elements (i.e.\ morphisms 
$1 \to \Sigma$) $\top$ and $\bot$. The set of global elements of $D$ (which
we also denote by $D$) can be endowed with the structure of a pca as follows:
for $t,s \in D$ we define $ts \in D$ as $(ts)(\vec{r}) = t(s.\vec{r})$. For
the set $\Lambda$ of terms we take $D$ and for the set $\Pi$ of stacks we take
$D^\omega$. The push operation sends $t \in \Lambda$ and $\vec{s} \in \Pi$ to
$t.\vec{s}$, the stream with head $t$ and tail $\vec{s}$. For every 
$\vec{s} \in \Pi$ let $\kk_{\vec{s}} \in \Lambda$ be defined as 
$\kk_{\vec{s}}(t.\vec{r}) = t(\vec{s})$. The control operator $\Cc$ is given
by $\Cc(t.\vec{s}) = t(\kk_{\vec{s}}.\vec{s})$. A natural choice for the pole 
$\Perp$ is $\{ \langle t , \vec{s} \rangle \mid t(\vec{s}) = \top \}$.

But on this level of generality we do not know how to choose a set $\PL$ of 
``proof-like terms''. However, in case  $D$ is the bifree solution of 
$D \cong \Sigma^{D^\omega}$ in some category of domains like 
\begin{enumerate}
\item[1)] cpo's with bottom and Scott continuous functions 
\item[2)] coherence spaces and stable (continuous) maps 
\item[3)] observably sequential algorithms as in \cite{ccf94}
\end{enumerate}
by a theorem of A.~Pitts (see \cite{pitts96}) there exists a unique subset 
$P$ of $D$ such that $t \in P$ 
iff $t(\vec{s}) = \bot$ for all $\vec{s} \in P^\omega$. Such a $P$ qualifies
as a set $\PL$ of proof-like terms since $P$ is closed under application,
contains all elements definable in untyped $\lambda$-calculus and we also have
$\Cc \in P$.

For later use we remark that the identity map on $D$ is represented by 
$\ic \in P$ with $\ic(t.\vec{s}) = t(\vec{s})$.

\section{Some triposes induced by $(D,P)$}\label{sometrip}

Since $P$ is a subpca of the pca $D$ we may consider the 
\emph{relative realizability} topos $\E = \RT(D,P)$ 
induced by the tripos $\mathscr{P}$ over $\Se$ 
where for a set $I$ the fibre $\mathscr{P}^I$ is the preorder 
$\bigl(\Pow(D)^I,\vdash_I\bigr)$ with $\phi \vdash_I \psi$ iff
$\exists t \in P.\forall i \in I.\forall s \in \phi_i.\, ts \in \psi_i$
and for $u : J \to I$ reindexing along $u$ is given by precomposition 
with $u$. For the set $\Sigma_{\mathscr{P}}$ of propositions of $\mathscr{P}$ 
we may take $\Pow(D)$ and for the \emph{truth predicate} on 
$\Sigma_{\mathscr{P}}$ we may take $\id_{\Pow(D)}$. 

Notice that $\Sigma_{\mathscr{P}}$ contains an ``intermediate'' truth value 
$U = \{\top_D\}$ which is neither equivalent to 
$\bot_{\Sigma_{\mathscr{P}}} = \emptyset$ nor to $\top_{\Sigma_{\mathscr{P}}} = D$.
Moreover, in $\RT(D,P)$ the proposition $U = \{\top_D\}$ is equivalent to
$U \vee \neg U$ (since $\neg U = \emptyset$) but not to $D$. Thus 
$U \vee \neg U$ does not hold in $\RT(D,P)$ for which reason the topos 
$\RT(D,P)$ is not boolean. However, the truth value $U$ gives rise to 
the (Lawvere-Tierney) topology $j_U$ on $\Sigma_{\mathscr{P}} = \Pow(D)$ 
which is defined as $j_U(A) = (A \to U) \to U$  for $A \in \Pow(D)$.
We may form the full subtripos $\mathscr{P}_U$ of $\mathscr{P}$  consisting of 
$j_U$-closed predicates, i.e.\ $\phi \in \Pow(D)^I$ with
$j_U \circ \phi \vdash_I \phi$. Since $j_U = \neg_U \circ \neg_U$ 
with $\neg_U A = A \to U$ the fibres of $\mathscr{P}_U$ are all boolean. We 
write $\E_U = \RT(D,P)_U$ for the ensuing boolean subtopos of $\E = \RT(D,P)$.

As described in the previous section $P \subseteq D$ gives rise to a 
classical realizability structure with \emph{pole} 
$\Perp\; = \{ \pair{t}{\vec{s}} \mid t(\vec{s}) = \top\}$. 
We write $\E_\Perp = \RT(D,P)_\Perp$ or rather simply $\Kcal$ for the ensuing 
classical realizability topos which is induced by the full subtripos 
$\mathscr{P}_\Kcal$ of $\mathscr{P}$ consisting of those predicates 
$\phi \in \Pow(D)^I$ which factor through 
$\Sigma_\Kcal = \{ A \in \Pow(D) \mid A^{\Perp\Perp} = A \}$. 
We show now that 

\begin{lem}\label{equivlm1}
$\mathscr{P}_\Kcal$ is equivalent to $\mathscr{P}_U$.
\end{lem}
\proof
First recall that on $\Pow(D)$ implication is given by
$A \to B = \{ t \in D \mid \forall s \in A.\, ts \in B \} 
 = \{ t \in D \mid \forall s \in A.\,\lambda \vec{r}. t(s.\vec{r}) \in B \}$
from which it follows that $\Sigma_\Kcal$ is an exponential ideal in $\Pow(D)$, 
i.e.\ $A \to B$ is in $\Sigma_\Kcal$ whenever $B$ is in $\Sigma_\Kcal$. 
Since $U \in \Sigma_\Kcal$ the map $j_U$ sends $\Pow(D)$ to $\Sigma_\Kcal$. 
Thus, postcomposition with $j_U$ gives rise to a tripos morphism from 
$\mathscr{P}$ to $\mathscr{P}_\Kcal$ left adjoint to the inclusion of tripos 
$\mathscr{P}_\Kcal$ into the tripos $\mathscr{P}$ (as induced by 
$\Sigma_\Kcal \subseteq \Pow(D)$). 
Since $A \to j_U(A)$ is uniformly realized by $\eta = \lambda x.\lambda p. px
\in P$ and for $A \in \Sigma_\Kcal$ the implication $j_U(A) \to A$ is uniformly 
realized by $\Cc \in P$ the adjunction above between $\mathscr{P}$ and 
$\mathscr{P}_U$ restricts to an equivalence between $\mathscr{P}_U$ and 
$\mathscr{P}_\Kcal$.\footnote{\emph{Question} We know that 
$j_U(A) \to A^{\Perp\Perp}$ is realized by $\Cc$ uniformly in $A \in \Pow(D)$. 
But is the reverse implication also realizable uniformly in $A$?}
\qed

\emph{Thus $\Kcal = \RT(D,P)_\Perp$ and $\RT(D,P)_U$ are equivalent boolean 
subtoposes of the relative realizability topos $\E = \RT(D,P)$ which itself
is not boolean.} 
We write $i : \Kcal \inc \E$ for the corresponding injective geometric 
morphism. Its inverse image part $i^* : \E \to \Kcal$ (sheafification) is 
given by postcomposition with $j_U$. Its (right adjoint) direct image part 
$i_* : \Kcal \to \E$ is nontrivial. As described in \cite{jaapbook} it sends 
an object $X$ in $\Kcal$ to $S(X)$, the object of singleton predicates on $X$ 
in $\Kcal$ considered as an object of $\E$.

Notice, however, that interpretation in $\E_U$ is not the same as (a kind of)
negative translation (with $U$ for falsity) followed by interpretation in $\E$
since this leaves interpretations of terms unchanged. But the finite type
hierarchies over natural numbers are quite different in $\E$ and in $\E_U$
from type level 2 onwards because of the presence of discontinuous 
functionals in $\E_U$ as opposed to $\E$ where all functionals over natural 
numbs are continuous.

For convenience and later use we explicitate a bit the logical structure 
of the triposes introduced above. 

For $A,B \in \Pow(D)$ implication in $\mathscr{P}$ is given by 
$A \to B = \{ t \in D \mid \forall s \in A.\, ts \in B\}$. 
Since the local operator $j_U$ commutes with this implication it also works
for $\mathscr{P}_U$. Looking a bit closer one sees that this holds also for
$\mathscr{P}_\Kcal$ since if $A$ and $B$ are biorthogonally closed then 
$A \to B = \{ t \in D \mid \forall s \in A.\forall \vec{r} \in B^\Perp.
 t(s.\vec{r}) = \top \} = 
\{ s.\vec{r} \mid s \in A , \vec{r} \in B^\Perp\}^\Perp$
and thus is biorthogonally closed, too.

For a set $I$ universal quantification $\forall_I$ along the terminal projection
$I \to 1$ is given by intersection, i.e.\ 
$\forall_I(\phi) = \bigcap\limits_{i \in I} \phi_i$. Since 
$\forall_I(\phi \to U) =  \bigl(\bigcup\limits_{i \in I} \phi_i\bigr) \to U$
it is immediate that $\forall_I$ restricts to $\mathscr{P}_U$. This applies
also to $\mathscr{P}_\Kcal$ since $\forall_I(\phi) = 
\bigcap\limits_{i \in I} \phi_i = \bigcap\limits_{i \in I} \phi_i^{\Perp\Perp} =
\bigl(\bigcup\limits_{i \in I} \phi_i^\Perp\bigr)^\Perp$ for which reason
$\forall_I(\phi)$ is biorthogonally closed. Universal quantification along
arbitrary maps $u : J \to I$ in $\Set$ is given by 
\[ \forall_u(\phi)_i = 
   \forall_J(\lambda j{\in}J. \mathit{leq}(u(j),i) \to \phi_j) \]
where $\mathit{leq}$ stands for Leibniz equality. 

Recall that Leibniz equality on set $I$ is defined as
\[ \mathit{leq}_I(i,j) = \bigcap\limits_{p \in \Sigma^I} p(i) \to p(j) \]
where $\Sigma$ refers to the $\Sigma$ of the respective tripos.
For the tripos $\mathscr{P}$ Leibniz equality on a set $I$ is given by 
$\mathit{leq}_I(i,j) = \{ \ic \mid i = j \}$. Obviously, the predicate
$\mathit{leq}_I$ is equivalent to the predicate $\eq_I$ defined as
$\eq_I(i,j) = \{ d \in D \mid i = j\}$. This observation is useful for
obtaining a simple description of equality predicates for the tripos 
$\mathscr{P}_U$ since they are of the form $j_U \circ \eq_I$.
Notice that $j_U(\emptyset) = (\emptyset \to U) \to U = 
D \to U = U = \{\top_D\}$ and $j_U(D) = (D \to U) \to U = U \to U 
= \{ d \in D \mid d \top_D = \top_D\} 
= \{ d \in D \mid \forall \vec{s} \in D^\omega.\, d(\top_D.\vec{s}) = \top \}
= \{ \top_D \} \cup {\uparrow}\bar{0}$ where $\bar{0}$ is the least element
of $D$ sending $\top_D.\bot_D^\infty$ to $\top$.
Thus, for $\mathscr{P}_U$ equality on $I$ is given by
$\eq_I(i,j) = \{\top_D\} \cup {\uparrow}\{ \bar{0} \mid i = j \}$.
A different but equivalent implementation of equality on $I$ for $\mathscr{P}_U$
is given by $\eq_I(i,j) = \{ \top_D \} \cup \{ \ic \mid i=j \}$ since there is
a least $r \in P$ with $r\bot_D = \bot_D$, $r\top_D = \top_D$ and $rd = \ic$
for $d \sqsupseteq \bar{0}$.

Since $\mathscr{P}_\Kcal$ is equivalent to its subtripos $\mathscr{P}_U$ the 
above considerations apply to $\mathscr{P}_\Kcal$ as well. 
 
\section{Nothing new in case of Scott domains}

In a talk in Chambery in June 2012 \cite{kriv2012} Krivine has shown 
that a classical realizability model is a forcing model iff it validates 
the sentence\footnote{in our terminology this means that
$\forall x{:}2. (eq_2(x,0) \vee \eq_2(x,1))$
holds in the tripos $\mathscr{P}_\Kcal$}
$\forall x^{\gimel 2} (x \neq 0, x \neq 1 \to \bot)$, i.e.\ iff there exists 
a proof-like term realizing $|\top,\bot \to \bot| \cap |\bot,\top \to \bot|$. 
He has shown that from such a realizer one can construct a proof-like term $\Phi$ 
such that $\Phi \in |A|$ whenever $|A|$ contains some proof-like term.

This applies in particular to the realizability structures 
as described in section~\ref{rsd} where $|\top| = D$ and $|\bot| = \{\top_D\}$.
Obviously, in this case $t \in P$ realizes 
$|\top,\bot \to \bot| \cap |\bot,\top \to \bot|$
iff $t \top_D s = \top_D = t s \top_D$ for all $s \in D$. But since $t \in P$
entails $t \bot_D \bot_D \neq \top_D$ this would give rise to a morphism 
$\vee : \Sigma \times \Sigma \to \Sigma$ with $u \vee v = \bot$ iff $u = v = \bot$
which does not exist in stable domain theory. However, in Scott domains such a
morphism does exist (``parallel or'') and allows one to construct an element of $P$
realizing $|\top,\bot \to \bot| \cap |\bot,\top \to \bot|$.
Moreover, in the case of Scott domains the classical realizability model
induced by $D$ and $P$ is not only a forcing model but it is actually 
equivalent to the ``ground model'' $\Set$ as we show next.

Since $P$ is Scott closed and closed under binary suprema  
it contains a greatest element $\Phi = \bigsqcup P$. Obviously, we have 
$\Phi(\vec{s}) = \bot$ iff $\vec{s} \in P^\omega$. Thus, a proposition $A$ holds 
in the ensuing realizability model (i.e.\ $|A| \cap P \neq \emptyset$) iff 
$\Phi \in |A|$ (since $|A| = ||A||^\Perp$ is upward closed).
Now for propositions $A$ and $B$ we have 
\begin{tabbing}
\qquad \= $\Phi \in |A \to B|$  \quad iff\\ 
\> $\forall t \in |A| \forall \vec{s} \in ||B||\; \Phi(t.\vec{s}) = \top$ \quad iff\\
\> $\forall t \in |A| \forall \vec{s} \in ||B|| \; 
   t \not\in P \vee \vec{s} \not\in P^\omega$ \quad iff\\
\> $\forall t \in |A|  \forall \vec{s} \in ||B|| \; 
    t \in P \Rightarrow \vec{s} \not\in P^\omega$ \quad iff\\
\> $\forall t \in |A| \; \bigl(t \in P \Rightarrow 
    \forall \vec{s} \in ||B|| \; \vec{s} \not\in P^\omega \bigr)$ \quad iff\\
\>  $\forall t \in |A| \; \bigl(t \in P \Rightarrow \Phi \in |B|\bigr)$ \quad iff\\
\> $\bigl(\exists t \in P \; t \in |A|\bigr) \Rightarrow \Phi \in |B|$ \quad iff\\
\> $\Phi \in |A| \Rightarrow \Phi \in |B|$
\end{tabbing} 
i.e.\ $A \to B$ holds iff from validity of $A$ follows validity of $B$.
Thus the ensuing classical realizability model is a 2-valued forcing model, 
i.e.\ coincides with the ground model $\Set$.

The situation changes dramatically if one solves the domain equation for $D$ 
in a category not admitting $\vee : \Sigma \times \Sigma \to \Sigma$ as e.g.\ 
the category $\mathbf{Coh}$ of coherence spaces and stable maps 
(see \cite{prftypbook}), the category $\mathbf{OSA}$ of observably sequential 
algorithms (see \cite{ccf94}) or a category of HON games and innocent algorithms. 
Let us look more closely at
the example of $D = \Sigma^{D^\omega}$ in $\mathbf{Coh}$ in which 
$||\bot|| = D^\omega$ and $||\top|| = \emptyset$ and accordingly 
$|\bot| = \{\top_D\}$ and $|\top| = D$.
Now if $f \in |\top,\bot \to \bot| \cap |\bot,\top \to \bot|$ 
then $f \top_D \bot_D = \top_D = f \bot_D \top_D$ and thus, by stability 
of $f$, also $f \bot_D \bot_D = \top_D$ from which it follows that 
$f \not\in P$. Thus the ensuing classical realizability model cannot be a 
forcing model (sheaves over a complete Boolean algebra) and, accordingly, 
is in particular different from the ground model $\Set$.

\section{Bifree Solution of $D = \Sigma^{D^\omega}$ in $\mathbf{Coh}$}

Let $V$ be the least set with $V = \Pfin(\omega{\times}V)$. 
If $\alpha \in V$ and $n\in\omega$ we write $\alpha_n$ for the set 
$\{ \beta \in V \mid \pair{n}{\beta} \in \alpha \}$.
By recursion on $n\in\omega$ we will define a sequence of coherence spaces
$D_n = \bigl(|D_n|,\coh_n\bigr)$ with $|D_n| \subseteq V$ such that 
\begin{enumerate}
\item[(1)] $|D_n| \subseteq |D_{n+1}|$ 
\item[(2)] for $\alpha,\beta\in|D_n|$ we have 
           $\alpha \coh_n \beta$ iff $\alpha \coh_{n+1} \beta$
\item[(3)] $\alpha \incoh_n \beta$ iff $\alpha \cup \beta \in D_n$
\item[(4)] $\alpha \coh_n \beta$ iff 
           $\alpha \cup \beta \in D_n$ implies $\alpha = \beta$. 
\end{enumerate}

For getting the construction of the $D_n$ right it is usful to recall that
coherence spaces and linear continuous maps between them give rise\footnote{Actually, this model was the \emph{source} of linear logic! But classical realizability does not make any use of the fact that $D$ lives within the Kleisli category
of the ``archetypical'' coherence space model for linear logic. Actually, we rather think of classical realizability as a non-linear version of Girard's \emph{Ludics}.} to a model of linear logic 
and that $\Sigma^{D^\omega} =\;!(D^\omega) \lolli \bot = (!(D^\omega))^\perp$. 
We put $|D_0| = \emptyset$, i.e.\ $D_0$ is the terminal object in 
$\mathbf{Coh}$. Notice that (3) and (4) vacuously hold for $D_0$.
For the induction step we put 
$D_{n+1} = \Sigma^{D_n^\omega} = \bigl(!(D_n^\omega)\bigr)^\perp$ as suggested by
$\Sigma^{D^\omega} =\;!(D^\omega) \lolli \bot$.
Thus, the web $|D_{n+1}|$ of $D_{n+1}$ consists of all $\alpha \in V$ 
such that for all $k\in\omega$ it holds that $\alpha_k \in D_n$, i.e.\ 
$\beta \coh_n \gamma$ for all $\beta,\gamma \in \alpha_k$, since $|D_{n+1}|$ 
is the web of $!(D_n^\omega)$ and for this coherence space we have 
$\alpha \coh \beta$ iff $\alpha \cup \beta \in |D_{n+1}|$. Thus, for defining
its orthogonal $D_{n+1}$ we put $\alpha \coh_{n+1} \beta$ iff 
$\alpha \cup \beta \in |D_{n+1}|$ implies $\alpha = \beta$. Conditions
(3) and (4) hold for $D_{n+1}$ by construction since they hold for $D_n$ by
induction hypothesis. We write $D$ for the coherence space where 
$|D| = \bigcup\limits_{n\in\omega} |D_n|$ and 
$\coh$ is the union of the $\coh_n$. 

Actually, one can avoid any explicit reference to the levels $D_n$ and 
inductively define $|D|$ as the least subset of $V$ with $\alpha \in |D|$ 
whenever $\forall n{\in}\omega.\,\alpha_n \subseteq |D| \wedge 
          \forall \beta,\gamma \in \alpha_n.\,\beta\coh\gamma$ 
where $\beta\coh\gamma$ stands for 
$\beta\cup\gamma \in |D| \Rightarrow \beta = \gamma$. 
Notice that $|D|$ is closed under subsets and we have 
$\alpha \incoh \beta$ iff $\alpha\cup\beta \in |D|$.

Now we describe the realizability structure arising from $D$.
The elements of $\Lambda_D = D$ are those $t \in \Pow(|D|)$ such that
$\forall \alpha,\beta \in t.\,\alpha \coh \beta$, 
i.e.\ \emph{antichains} in the poset $(|D|,\subseteq)$. 
The evaluation map $D \times D^\omega \to \Sigma$ is defined as follows: 
for $t \in D$ and $\vec{s} \in D^\omega$ we have $t(\vec{s}) = \top$ 
(notation $t \star \vec{s} \in\,\Perp$) iff 
$\exists \alpha \in t.\forall n \in\omega.\,\alpha_n \subseteq s_n$.
With an $\vec{s} \in D^\omega$ one may associate the set 
$I_{\vec{s}} = \{ \alpha \in |D| \mid \{\alpha\}(\vec{s}) = \top \} = 
 \{ \alpha \in |D| \mid \forall n\in\omega.\, \alpha_n\subseteq s_n\}$.
Sets of this form can be characterized as downward closed ideals in $|D|$,
i.e.\ subsets of $D$ which are closed under subsets and finite unions. Any 
such ideal $I$ is equal to $I_{\vec{s}}$ for a unique $\vec{s} \in D^\omega$ 
which is given by $s_n = \bigcup\limits_{\alpha\in I} \alpha_n$. Writing
$\Pi_D$ for the set of downward closed ideals in $(|D|,\subseteq)$ 
for $t \in \Lambda_D$ and $\pi \in \Pi_D$ we have 
$t \star \pi \in\;\Perp$ iff $t \cap \pi \neq \emptyset$.

For exhibiting in a concrete way the remaining operations of the realizability 
structure induced by $D = \Sigma^{D^\omega}$ we have to introduce some notation. 
For a finite $a \in D$ and $\alpha \in |D|$ we write $a.\alpha$ for 
$(\{0\} \times a) \cup \{ \pair{n+1}{\beta} \mid \pair{n}{\beta} \in\alpha \}$.
For $t \in \Lambda_D$ and $\pi \in \Pi_D$ let 
$t.\pi = \{ a.\alpha \mid  a \subseteq_\fin t , \alpha \in \pi \}$. 
For $t,s \in \Lambda_D$ let $ts = \{ \alpha \in |D| \mid 
\exists a \subseteq_\fin s.\; a.\alpha \in t \}$. For $f \in \mathbf{Coh}(D,D)$
let $\mathsf{fun}(f) = \{ a.\alpha \mid (a,\alpha) \in \mathsf{tr}(f)\}$ where
$\mathsf{tr}(f)$ is the trace of $f$, i.e.\ the set of all pairs $(a,\alpha)$ 
s.t.\ $a \in D$ is finite and $\alpha \in f(a)$ and for all $b \subseteq a$ 
from $\alpha \in f(b)$ it follows that $a = b$. Using $\mathsf{fun}$ we define
$\lambda x .t = \mathsf{fun}(a \mapsto t[a/x])$. For $\pi\in\Pi_D$ we put
$\kk_\pi = \{ \widehat{\alpha} \mid \alpha \in \pi \}$ 
where $\widehat{\alpha} = \{\pair{0}{\alpha}\}$ for $\alpha \in |D|$.
We define $\Cc = \bigl\{ \{\{\widehat{\alpha}_1,\dots,\widehat{\alpha}_k\}.\alpha\}.(\alpha{\cup}\alpha_1{\cup}\dots{\cup}\alpha_k) \mid 
\alpha{\cup}\alpha_1{\cup}\dots{\cup}\alpha_k \in |D| \bigr\}$. 

Finally, we have to define which elements of $\Lambda_D$ we want to consider 
as proof-like objects. By recursion on $\alpha \in |D|$ we define 
$|\alpha| \in \{0,1\}$ as
$|\alpha| = 1$ iff $\exists n\in\omega.\exists\beta\in\alpha_n.\,|\beta| = 0$. 
Thus $|\alpha| = 1$ iff $\alpha$ does not raise any error itself. Accordingly,
we define the subset $P$ of \emph{proof-like} objects of $D$ as
$\{ a \in D \mid \forall \alpha \in a.\, |\alpha| = 1 \}$.

\subsection{Some useful retractions in $P$}

Let $h_n$ be the sequence of subidentical retractions of $D$ where 
$h_0 = \lambda f{:}D.\bot_D$ and $h_{n+1} = \Sigma^{h_n^\omega}$. Notice,
that all $h_n$ are elements of $P$ since $\id_D$ is in $P$ and 
$h_n \sqsubseteq \id_D$. Obviously, we have 
$\id_D = \bigsqcup\limits_{n\in\omega} h_n$ but the images of the $h_n$
typically contain elements which are not finite. Notice that the image
of $h_n$ is $D_n$ for all $n\in\omega$.

There is also a subidentical retraction $r_P \in P$ sending $a \in D$ to
$r_P(a) = \{ \alpha \in a \mid |\alpha| = 1 \}$. Obviously, the image of
$r_P$ is precisely $P$ and $r_P(a)$ is the greatest element of $P$ below $a$.

\subsection{$D$ is universal for countably based coherence spaces}

To give an impression of the complexity of $D$ we show that it contains every
countably based coherence space via a stably continuous embedding/projection pair
(see e.g.\ \cite{AspertiLongoBook}). First recall that in Th.~2.4.2.9 of 
\cite{AspertiLongoBook}) it has been shown that every coherence space $X$ 
with countable web can be embedded into $!\mathbb{T}^\omega$ via a stably continuous 
embedding/projection pair (where $\mathbb{T}$ is the coherence space whose web 
consists of two incoherent tokens thought of as boolean values). Thus, the
coherence space $X^\perp$ can be embedded into 
$(!\mathbb{T}^\omega)^\perp = \Sigma^{\mathbb{T}^\omega}$. Accordingly, all 
coherence spaces with countable web can be embedded into $\Sigma^{\mathbb{T}^\omega}$. 
Since $\mathbb{T}$ can be embedded into $D$ the coherence space 
$\Sigma^{\mathbb{T}^\omega}$ can be embedded into $\Sigma^{D^\omega}$ and thus into $D$.

\subsection{Antichains in Coherence Spaces}

Let $X$ be a coherence space. An \emph{antichain} in $X$ is a subset $A$ 
of $X$ such that $a,b \in A$ are equal whenever they are coherent 
(i.e.\ $a \cup b \in X$). We may order antichains in $X$ ``\`a la Smyth'' 
as follows
\[  A \leq_S B \quad\mbox{iff}\quad 
    \forall y \in B. \exists x \in A.\; x\sqsubseteq y \]
i.e.\ $A \leq_S B$ iff ${\uparrow}A \supseteq {\uparrow}B$. This suggest 
to consider antichains as upward closed subsets $C$ of $X$ such that 
for the set $\min(C)$ of minimal elements of $C$ it holds that
\begin{enumerate}
\item[(1)] $C \subseteq {\uparrow}\min(C)$ and
\item[(2)] coherent elements of $\min(C)$ are equal.
\end{enumerate}
Under this view antichains may be considered as disjoint unions of 
\emph{cones}, i.e.\ sets of the form ${\uparrow}x$ for some $x \in X$.
We write $\mathscr{A}(X)$ for the set of antichains considered as upward
closed subsets of $X$ satisfying conditions (1) and (2) and consider it
partially ordered by reverse subset inclusion. One can show that 
\begin{thm}
$\mathscr{A}(X)$ is a complete lattice when ordered by $\supseteq$.
\end{thm}
\proof
Let $(C_i)_{i \in I}$ be a family of antichains in $X$. We show that its
intersection $D := \bigcap\limits_{i \in I} C_i$ is again an antichain 
from which it is immediate that $D$ is the supremum of the $C_i$ 
w.r.t\ $\supseteq$.

Obviously, the set $D$ is upwards closed. For $x \in D$ and $i \in I$
let $x_i$ be the unique element of $\min(C_i)$ with $x_i \sqsubseteq x$.
Since $(x_i)_{i \in I}$ is bounded by $x$  its supremum $m(x)$ exists. 
It is easy to see that $x \sqsupseteq m(x) \in \min(D)$. Thus $D$ validates
condition (1). For showing condition (2) suppose $x,y \in \min(D)$ have an
upper bound. Then for all $i \in I$ we have $x_i \coh y_i$ and thus 
$x_i = y_i$ from which it follows that $m(x) = m(y)$ and thus $x = y$ 
as desired.
\qed

An important class of antichains in $X$ are those of the form $p^{-1}(\top)$
for some $p \in \Sigma^X$. Via trace they correspond to those 
$U \in \mathscr{A}(X)$ for which all elements of $\min(U)$ are compact
elements of $X$. We write $\mathscr{A}_0(X)$ for this class of antichains 
in $X$. 
For every $C \subseteq X$ we may consider the antichain
\[  \overline{C} = \bigcap \{ U \in \mathscr{A}_0(X) \mid U \supseteq C \} \]
which, obviously, contains $C$ as a subset. It is easy to see that
$C \mapsto \overline{C}$ is a closure operator on $\Pow(X)$ since $\overline{C}$
is the intersection of all stably open subsets of $X$ which contain $C$ 
as a subset.\footnote{If $X$ is a Scott domain then intersections of open subsets
of $X$ are just upward closed subsets of $X$. Alas, such an easy characterizations
is not available for intersections of stably open subsets of a coherence space $X$.}

For $X = D^\omega$ and $C \subseteq X$ we have 
$C^\Perp = \{ t \in \Sigma^X \mid C \subseteq t^{-1}(\top) \}$ 
and thus $\overline{C} = C^{\Perp\Perp}$. Notice that the minimal elements of $C^\Perp$ 
w.r.t.\ the stable order are those $t \in D = \Sigma^{D^\omega}$ for which every 
element of $\trace(t)$ is below some element of $C$. For $t \in C^\Perp$ the unique 
minimal element $m(t)$ in $C^\Perp$ below $t$ is characterized as follows: 
$e \in \trace(m(t))$ iff $e \in \trace(t)$ and $e \sqsubseteq x$ for some $x \in C$. 

Since the infimum operation $\sqcap : \Sigma \times \Sigma \to \Sigma$ is stably
continuous for $t_1,t_2 \in D$ we have $t_1 \sqcap t_2 \in D$. Obviously, 
we have $(t_1 \sqcap t_2)^{-1}(\top) = t_1^{-1}(\top) \cap t_2^{-1}(\top)$. Thus $C^\Perp$
is not only an antichain in $D$ but it is also closed under $\sqcap$ and
contains $\lambda x{:}X.\top$ as an element. It is an interesting but difficult 
problem to characterize those antichains $A$ in $D$ which are of the form
$C^\Perp$ for some $C \subseteq X$. Well, it are those $A \subseteq D$ for
which $A = A^{\Perp\Perp}$. But is there a more elementary combinatorial 
characterization of biorthogonally closed subsets of $D$? Such a characterization
might be helpful for answering the question whether for any biorthogonally closed
subset $A$ of $D$ either $A$ or its negation $\neg_UA$ is inhabited by an element  
of $P$, i.e.\ whether $\Kcal$ is 2-valued.

\section{Exploring the structure of $\E$ and $\Kcal$}

We have seen that $\Kcal$ is equivalent to $\Set$ when constructed from the 
bifree solution of $D = \Sigma^{D^\omega}$ in Scott domains. But something new
arises when we start from the solution of this domain equation 
in $\mathbf{Coh}$. We start now exploring this new territory. Some attention
will also be payed to the inituitionistic variant $\E$ in which computation
is much easier than in its full subcategory $\Kcal$ of $j_U$-sheaves.

For every $n\in\N$ let $\bar{n}$ be the unique element of $D$ with
$\bar{n}(\vec{s}) = \top$ iff $s_n = \top_D$. A ``hardwired'' version of this
is $\bar{n} = \{\nu_n\}$ with $\nu_n = \{\pair{n}{\emptyset}\} \in |D|$.
From this it is obvious that the $\bar{n}$ are atoms of $D$ and pairwise 
incoherent, i.e.\ $\nu_n \sincoh \nu_m$ iff $n \neq m$. Obviously, we have
$\bar{n} \in P$ since $|\nu_n| = 1$.

In $\E = \RT(D,P)$ a natural numbers object is given by the assembly $N_\E$
with underlying set $\N$ and $||n||_{N_\E} = \{\bar{n}\}$. Similarly, the 
object $2_\E$ in $\E$ is given by the assembly with underlying set 
$2  = \{0,1\}$ and $||k||_{2_\E} = \{\bar{k}\}$. 
The object $\Delta_\E(2)$ of $\E$ is given by the assembly with underlying 
set $2$ and $||k||_{\Delta_\E(2)} = D$.

The corresponding objects $N_\Kcal$, $2_\Kcal$ and $\Delta_\Kcal(2)$ in $\Kcal$ 
are obtained from $N_\E$, $2_\E$ and $\Delta_\E(2)$ in $\E$ by 
\emph{sheafification} (denoted as $i^*$), i.e.\ by postcomposing the 
respective equality predicates with $j_U$. But since $j_U$ is a bit complex 
we are looking for somewhat simpler isomorphic copies of these objects 
in $\Kcal$.

Since $j_U(\emptyset) = \{\top_D\}$ and 
$j_U(D) = \{\top_D\} \cup {\uparrow}\{\bar{0}\}$ for every set $I$ 
the object $\Delta_\Kcal(I)$ of $\Kcal$ has underlying set $I$ and 
equality predicate $\sem{i \sim_{\Delta_\Kcal(I)} j} = 
\{\top_D\} \cup  {\uparrow}\{ \bar{0} \mid i=j \}$.

Next we determine $i^*N_\E$, the natural numbers object of $\Kcal$ obtained
by sheafifying the natural numbers object $N_\E$ of $\E$. The underlying
set of $i^*N_\E$ is $\N$ and its equality predicate is given by 
$\sem{n \sim_{i^*N_\E} m} = j_U\bigl(\sem{ n \sim_{N_\E} m}\bigr)$.
The following lemma exhibits an object $N_\Kcal$ which in $\Kcal$ is isomorphic
to  $i^*N_\E$ but simpler to describe and simpler to manipulate.

\begin{lem}\label{NKlm}
Let $N_\Kcal$ be the object of $\Kcal$ with underlying set $\N$ and equality
predicate $\sem{ n \sim_{N_\Kcal} m}  = 
\{\top_D\} \cup {\uparrow}\{ \bar{n} \mid n = m\}$. In $\Kcal$ the object
$N_\Kcal$ is isomorphic to $i^*N_\E$ and thus a natural numbers object 
in $\Kcal$
\end{lem}
\proof
For showing the desired isomorphism it suffices to exhibit elements of $P$
realizing the logical equivalence of $\sem{ n \sim_{N_\Kcal} m}$ and
$j_U\bigl(\sem{ n \sim_{N_\E} m}\bigr)$ uniformly in $n$ and $m$.

First notice that 
$\sem{ n \sim_{N_\Kcal} n}  = 
\{\top_D\} \cup {\uparrow}\bar{n} = 
\{ \vec{s} \in D^\omega \mid s_n = \top_D\}^\Perp \in \Sigma_\Kcal$.
Next we determine $j_U(\{\bar{n}\}) = (\{\bar{n}\} \to U) \to U$ for $n \in\N$. 
Observe that $\{\bar{n}\} \to U = \{ d \in D \mid d\bar{n} = \top_D\} =
\{ d \in D \mid \forall \vec{s} \in D^\omega.\, d(\bar{n}.\vec{s}) = \top\}$.
Thus, we have $j_U(\{\bar{n}\}) = \{ d \in D \mid \forall d^\prime \in D.\, 
d^\prime\bar{n} = \top_D \Rightarrow dd^\prime = \top_D\}$.

First we show that $\Cc \in P$ realizes the implication 
$j_U\bigl(\sem{ n \sim_{N_\E} m}\bigr) \to \sem{ n \sim_{N_\Kcal} m}$
uniformly in $n$ and $m$. If $n \neq m$ then 
$\sem{n \sim_{N_\Kcal} m} = \{\top_D\} = j_U(\emptyset) =
j_U\bigl(\sem{ n \sim_{N_\E} m}\bigr)$ and the claim follows 
since $\Cc\top_D = \top_D$. Thus it suffice to show that $\Cc$ realizes
$j_U\bigl(\sem{ n \sim_{N_\E} n}\bigr) \to \sem{ n \sim_{N_\Kcal} n}$ for all $n$.
For this purpose suppose $t \in j_U(\{\bar{n}\})$ and 
$\vec{s} \in D^\omega$ with $s_n = \top_D$. 
Then $\kk(\vec{s}) \in \{\bar{n}\} \to U$ since for $\vec{r} \in D^\omega$
we have $\kk(\vec{s})(\bar{n}.\vec{r}) = \bar{n}(\vec{s}) = \top$.
Thus, we have $\Cc(t.\vec{s}) = t(\kk(\vec{s}).\vec{s}) = \top$ as desired
since $t \in j_U(\{\bar{n}\})$ and $\kk(\vec{s}) \in \{\bar{n}\} \to U$. 

There is an $e \in P$ with $e \top_D = \top_D$ and 
$e\bar{n}d = d\bar{n}$ for all $n\in\N$ and $d \in D$. 
Obviously, such an $e$ realizes $\{\top_D\} \to j_U(\emptyset)$. 
Moreover, for $n\in\N$ we have 
$e\bar{n} \in j_U(\{\bar{n}\})$ since if $d\bar{n} = \top_D$ then also
$e\bar{n}d = d\bar{n} = \top_D$. Thus, since $j_U(\{\bar{n}\})$ is upward 
closed for every $d \sqsupseteq \bar{n}$ we have $e\bar{n} \sqsubseteq ed \in
j_U(\{\bar{n}\})$. Thus $e$ realizes $\{\top_D\} \cup {\uparrow}\bar{n} \to
j_U(\{\bar{n}\})$. Thus, we have shown that $e$ realizes 
$\sem{ n \sim_{N_\Kcal} m} \to j_U\bigl(\sem{ n \sim_{N_\E} m}\bigr)$
uniformly in $n$ and $m$.  
\qed

Similarly, one shows that in $\Kcal$ the object $i^*2_\E$ is isomorphic 
to the object $2_\Kcal$ with underlying set $2 = \{0,1\}$ and equality predicate
$\sem{i \sim_{2_\Kcal} j} = 
 \{\top^D\} \cup {\uparrow}\{\overline{i} \mid  i = j\}$.
Since $\Kcal$ is a boolean topos the truth value object $\Omega_\Kcal$ is
known to be isomorphic to $2_\Kcal$. We do not know whether the object 
$2_\Kcal$ has precisely two global elements, i.e.\ whether the topos $\Kcal$ 
is 2-valued.\footnote{But it can be shown that for \emph{countable} 
$A \subseteq D$ either $A^{\Perp\Perp}$ or its negation are true in $\Kcal$.}

Since $\Kcal$ is a subtopos of $\E$ arising from the Lawvere-Tierney topology
$j_U$ on $\E$ there is an induced injective geometric morphism $i : \Kcal \inc \E$
whose inverse image part $i^* : \E \to \Kcal$ we have already described. It is
fairly simple since it is given by postcomposition with $j_U$. However, its 
right adjoint $i_*$, the direct image part of $i$, though full and faithful is 
not simply inclusion in the naive sense. As described e.g.\ in \cite{jaapbook} 
it sends an object $X$ of $\Kcal$ to the object $i_*X$ of $\E$ which is the 
object $S(X)$ of `singleton predicates' on $X$ in $\Kcal$ considered as an 
object of $\E$. The underlying set of $S(X)$ is the set of all functions from 
$|X|$ to $\Sigma_\Kcal$ where $|X|$ is the underlying set of $X$. The existence 
predicate $E_{S(X)}$ on $\Sigma_\Kcal^{|X|}$ is given by
\[ E_{S(X)}(A) = \sem{\mathit{Pred}_X(A) \wedge \exists x{:}|X|.\forall y{:}|X|.
                     A(y) \leftrightarrow x \sim_X y} \]
where
\[ \mathit{Pred}_X(A) = \sem{\forall x{:}|X|. A(x) \to x \sim_X x \wedge 
                             (\forall y{:}|X|. x \sim_X y \to A(y))} \]
and the equality predicate for $S(X)$ is given by
\[ \sem{A \sim_{S(X)} B} = \sem{E_{S(X)}(A) \wedge 
                               \forall x{:}|X|. A(x) \leftrightarrow B(x)} \]  
which finishes the description of $S(X)$. For the morphism part of $S$ suppose
$F : |X| \times |Y| \to \Sigma_\Kcal$ represents a morphism from $X$ to $Y$. Then
the corresponding morphisms from $S(X)$ to $S(Y)$ is given by the 
$\Sigma_\Kcal$-valued predicate $S(F) : \Sigma_\Kcal^{|X|} \times \Sigma_\Kcal^{|Y|}
\to \Sigma_\Kcal$ defined as
\[ S(F)(A,B) = \sem{E_{S(X)}(A) \wedge E_{S(Y)}(B) \wedge \forall x{:}|X|,y{:}|Y|.
               F(x,y) \leftrightarrow (A(x) \wedge B(y))} \]
for $A \in \Sigma_\Kcal^{|X|}$ and $B \in \Sigma_\Kcal^{|Y|}$.
Thus, though the inclusion of $\Kcal$ into $\E$ via $i_*$ preserves 
exponentials due to the complicated nature of $i_*$ there is not much gain when
computing the exponentials in the relative realizability topos $\E$. 

Generally, since classical realizability toposes are boolean $\Omega$ is 
isomorphic to $2$. Thus, since $2$ is a subobject of $N$ the exponential $N^N$ 
contains $2^N \cong \Pow(N)$ as a subobject which explains why in general $N^N$ 
is so complicated in classical realizability toposes. Maybe this is the reason
why Krivine in his papers considers classical realizability models for classical
second order logic or the classical set theory $\mathsf{ZF}$ which are both based
on sets and not on functions. In both settings functions appear only as a derived
concept, namely as functional relations, i.e.\ particular sets.\footnote{Of course,
in case of second order logic he has to permit \emph{function constants} on the
underlying (countable) set of objects (usually identified with the set natural 
numbers).}

\smallskip
So far we do not know yet whether $\Kcal$ is actually different from a forcing model.
But it will follow from the results of the following subsection where we show that

\subsection*{$\Kcal$ is not even a Grothendieck topos}

Since there is no parallel-or in the realizability structure induced 
by $P \subseteq D$ it follows from Krivine's observation in \cite{kriv2012} 
that the object $\Delta_\Kcal(2)$ is not isomorphic to $2_\Kcal$. For this reason 
the tripos $\mathscr{P}_\Kcal$ does not arise from a complete boolean algebra. 
But from this it does not follow yet that $\Kcal$ is not equivalent to a forcing 
model, i.e.\ a localic boolean topos, since non-equivalent triposes might induce 
the same topos. But we will show now that $\Kcal$ is not even a Grothendieck topos 
and thus \emph{a fortiori} not a forcing model. 

For this purpose we will proceed in two steps. First in Lemma~\ref{subgrothtriv} 
we will show that every Grothendieck subtopos of $\Kcal$ is equivalent to $\Set$ 
and then in the subsequent Lemma~\ref{KnotSet} we will show that $\Kcal$ is not 
equivalent to $\Set$.
It is then an immediate consequence of these two lemmas that

\begin{thm}\label{KnotGroth}
$\Kcal$ is not a Grothendieck topos and thus, in particular, not a forcing model.
\end{thm}

The following considerations are necessary as preparation for the proofs of
Lemma \ref{subgrothtriv} and \ref{KnotSet}.

There is a geometric inclusion $\Pi_\E \dashv \Delta_\E : \Set \inc \E$ where
$\Pi_\E$ is given by $\E(U,-)$. The right adjoint $\Delta_\E$ sends set $I$ to 
the object $\Delta_\E(I) = (I,\eq_I)$ (see section~\ref{sometrip}) and 
$u : J \to I$ in $\Set$ to the morphism 
$\Delta_\E(u) : \Delta_\E(J) \to \Delta_\E(I)$ represented by the 
$\mathscr{P}$-predicate $\eq_I(u(j),i)$ on $J \times I$. Notice that $\Delta_\E$ 
factors through $\Asm(P,D)$, the category of assemblies in $\RT(D,P)$, since 
$\Delta_\E(I)$ is isomorphic to the assembly with underlying set $I$ and 
$||i|| = D$ for all $i \in I$. The restriction of the left adjoint $\Pi_\E$ 
to $\Asm(P,D)$ sends an assembly to its underlying set and a morphism to its 
underlying set-theoretic function. Notice that 
$\Pi_\E \dashv \Delta_\E : \Set \inc \E$ is the least non-trivial subtopos 
of $\E$ induced by the double negation topology on $\E$. 

We write $\bar{D}$ for the object of $\Asm(P,D)$ with underlying set $D$ and 
$||t||_{\bar{D}} = \{t\}$ for $t \in D$. Obviously, the counit 
$\eta_{\bar{D}} : \bar{D} \to \Delta_\E\Pi_E\bar{D}$ is monic. 
If $j : \F \inc \E$ is a nontrivial subtopos of $\E$ then the counit 
$\bar{D} \to j_*j^*\bar{D}$ of $j^* \dashv j_*$ at $\bar{D}$ factors 
along $\eta_D$ via a subobject $j_*j^*\bar{D} \mono \Delta_\E\Pi_E\bar{D}$ whose 
characteristic predicate $\chi_D$ is given by $\chi_D(t) = j_\F(\{t\})$ 
for $t \in D$ where $j_\F$ is the closure operator on $\E$ inducing the subtopos 
$\F$ of $\E$. Thus $j_*j^*\bar{D}$ is (isomorphic to) the assembly with underlying 
set $D$ and $||t||_{j_*j^*\bar{D}} = j_\F(\{t\})$. 

Now adapting an argument from \cite{ptj13} we show that

\begin{lem}\label{subgrothtriv}
Every nontrivial Grothendieck subtopos $\F$ of $\E$ is equivalent to $\Set$.
\end{lem}
\proof
Suppose $\F$ is a nontrivial Grothendieck subtopos of $\E$. We write 
$j : \F \inc \E$ for the corresponding inclusion. Since $\F$ is a Grothendieck 
topos it has arbitrary copowers.
We write $\Delta_\F(I)$ for the $I$-fold copower of $1_\F$, i.e.\ $\coprod_I 1_\F$.
Notice that $C = \Delta_\F(D)$ and $j_*j^*\bar{D}$ are both assemblies. Since 
$\Asm(P,D)$ is an exponential ideal in $\RT(D,P)$ the exponential 
$(j_*j^*\bar{D})^C$ is an assembly, too, and, moreover, (isomorphic to) the 
$D$-fold product of $j_*j^*\bar{D}$. The underlying set of $(j_*j^*\bar{D})^C$ 
may be identified with the set of all functions from $D$ to $D$ since
$\Pi_\E((j_*j^*\bar{D})^C) \cong \E(U,j_*j^*\bar{D})^C) \cong \E(C,(j_*j^*\bar{D})^U) 
\cong \E(U,j_*j^*\bar{D})^D \cong \Set(D,D)$. 

Let $j_\F$ be the closure operator on $\E$ giving rise to the subtopos $\F$ of $\E$.
The subobject $j_*j^*\bar{D} \mono \Delta_\E\Pi_\E \bar{D}$ is classified by the 
predicate $\chi_D(t) = j_\F(\{t\})$. Obviously, the Grothendieck topos $\F$ is 
equivalent to $\Set$ iff $\chi_D$ is constantly true, i.e.\ there is a $t \in P$ 
with $t \in \bigcap\limits_{t \in D} j_\F(\{t\})$.

For sake of contradiction suppose this were not the case. Then by axiom of choice 
on the meta-level there exists a (typically non-continuous) function $g : D \to D$ 
with $t \not\in j_\F(\{g(t)\})$. For $t \in D$ let $s_t \in P$ be some realizer 
for the projection $\pi_t : (j_*j^*\bar{D})^C \to j_*j^*\bar{D} : h \mapsto h(t)$. 
Let $f : D \to D : t \mapsto g(s_tt)$ for which it obviously holds that
$s_tt \not\in j_U(\{f(t)\})$ for $t \in D$. Since $\Pi_\E((j_*j^*\bar{D})^C) \cong 
\Set(D,D)$ there is a $t \in D$ realizing $f$ as an object of $(j_*j^*\bar{D})^C$.
But then $s_tt \in j_U(\{f(t)\})$ which is impossible.
\qed

Now for showing Theorem~\ref{KnotGroth} it remains to prove that

\begin{lem}\label{KnotSet}
The topos $\Kcal$ is not equivalent to $\Set$.
\end{lem}
\proof
For sake of contradiction suppose that $\Kcal$ is equivalent to $\Set$.
Then $i_*i^*\bar{D} \mono \Delta_\E\Pi_\E\bar{D}$ is an isomorphism. But then 
the predicate $\chi_D$ on $\Delta_\E\Pi_\E\bar{D}$ is constantly true, i.e.\ 
there is an $s \in P$ with $s \in j_U(\{t\})$ for all $t \in D$. But this
is impossible since already $j_U(\{\top_D\}) = (U \to U) \to U$ does not 
contain an element of $P$ (since such an element $t$ would map 
$\ic \in (U \to U) \cap P$ to an element $t\ic \in U \cap P = \emptyset$).
\qed

\section{$\Kcal$ is a model of full first order arithmetic}

Since $\Kcal$ hosts a natural numbers object $N_\Kcal$ it is most natural to ask
how much of first order arithmetic holds in $\Kcal$. First notice that \emph{all}
functions on $\N$ do exist as morphisms in $\Kcal$. An arbitrary set-theoretic
function $f : \N\to\N$ is represented as the morphism $f_\Kcal : N_\Kcal \to N_\Kcal$
as given by the $\mathscr{P}_\Kcal$-predicate $\sem{f(n) \sim_{N_\Kcal} m}$ 
on $\N \times \N$ because there exists $t_f \in P$ 
with $t_f \bar{n} = \overline{f(n)}$ for all $n\in\N$. Equality of natural
numbers will be interpreted as $\sem{\cdot \sim_{N_\Kcal} \cdot}$. 
Propositional logical connectives will be interpreted as usual 
(see \cite{kriv2009}) but notice that $|A{\to}B| = \{ t \in D \mid \forall s \in |A|.\, ts \in |B| \}$.\footnote{As 
in \cite{kriv2009} we write $|A|$ for the interpretation of formula $A$.} 
Universal quantification over $N_\Kcal$ is interpreted as 
$$|\forall x. A(x)| = 
  \bigcap\limits_{n\in\N} \sem{n \sim_{N_\Kcal} n} \to |A(n)|$$
which is coincidence with \cite{kriv2009} since the equivalence of 
$\sem{n \sim_{N_\Kcal} n}$ and 
\[ \bigcap\limits_{X \in \Sigma_\Kcal^\N} X(0) \to
\forall x (X(x) \to X(x{+}1)) \to X(n) \] 
can be realized by an element of $P$.
As usual existential quantification over $N_\Kcal$ is interpreted 
as its second order encoding, i.e.\ 
\[ |\exists x. A(x)| = \bigcap\limits_{X \in \Sigma_\Kcal} 
   \left(\bigcap\limits_{n\in\N} \bigl(\sem{n \sim_{N_\Kcal} n} \to |A(n)| \to X\bigr)\right) \to X \]
from which it follows that $\lambda f. f \bar{n} t \in |\exists x. A(x)|$ 
whenever $t \in |A(n)|$. Now we are ready to prove that

\begin{thm}
$\Kcal$ validates all true sentences of first order arithmetic.
\end{thm}
\proof
Since $\Kcal$ is boolean and classically every first order sentence is 
provably equivalent to a sentence in \emph{prenex form}, i.e.\ a prefix 
of quantifiers followed by an equation between arithmetic terms, it 
suffices to show that  all true arithmetic formulas in prenex form do hold 
in $\Kcal$. We proceed by structural induction on the structure of 
arithmetical sentences in prenex form.

If $e_1 = e_2$ is a true arithmetical equation where both sides have value $n\in\N$
then $e_1 = e_2$ is realized by $\bar{n} \in P$.

Suppose $\forall x. A(x)$ is a true arithmetical sentence in prenex form. Then 
for all $n\in\N$ the sentence $A(n)$ is true and in prenex form. Thus, by induction
hypothesis for every $n\in\N$ there is a $p_n \in P$ realizing $A(n)$. Then there
exists a $t \in P$ with $t\top_D = \top_D$ and $t\bar{n} = p_n$ for all $n\in\N$.
Obviously $t$ realizes $\forall x. A(x)$.

Suppose $\exists x. A(x)$ is a true arithmetical sentence in prenex form. Then 
for some $n\in\N$ the sentence $A(n)$ is true and in prenex form. By induction
hypothesis there is a $p \in P$ realizing $A(n)$ from which it follows that
$\lambda f. f \bar{n} p \in P$ realizes $\exists x. A(x)$.
\qed

Thus, w.r.t.\ first order arithmetic sentences one cannot distinguish $\Kcal$ from $\Set$.
But already at second order things get much more delicate since one does not even know 
whether every morphism $N_\Kcal \to N_\Kcal$ in $\Kcal$ is induced by a map $\N \to \N$
in $\Set$, i.e.\ whether for any functional relation $F$ from $N_\Kcal$ to $N_\Kcal$ 
there exists function $f : \N\to\N$ such that 
$\forall x,y{:}N_\Kcal.F(x,y) \leftrightarrow f(x) \sim_{n_\Kcal} y$ holds in $\Kcal$. One
easily sees that $f$ is uniquely determined by $F$ but the question rather is whether for
all $F$ such an $f$ exists. 

Actually, there is an even simpler question of this kind for which we do not know the answer 
so far, namely whether in $\Kcal$ the natural numbers object $N_\Kcal$ has only ``standard'' 
global elements. More explictly, this means whether for any morphism $a : 1_\Kcal \to N_\Kcal$ 
in $\Kcal$ there is an $n\in\N$ such that $a \sim_{N_\Kcal} n$ holds in $\Kcal$. The answer is
definitely negative for boolean valued models $\Sh(B)$ when $B$ is a complete boolean algebra
with more than $2$ elements. Since if $u \in B$ is different from $0_B$ and $1_B$ then so is 
$\neg u$ and one may cook up a ``mixed'' natural number $a$ which is $0$ on $u$ and $1$ 
on $\neg u$. We could come up with a similar ``nonstandard'' global element of $N_\Kcal$ 
if $\Omega_\Kcal = 2_\Kcal$ were not \emph{$2$-valued}, i.e.\ if their existed an 
$u : 1 \to 2_\Kcal$ for which $\Kcal$ validates neither $u \sim_{2_\Kcal} 0$ nor $u \sim_{2_\Kcal} 1$ though it certainly
validates the disjunction $u \sim_{2_\Kcal} 0 \vee u \sim_{2_\Kcal} 1$.

\section{The object $\Delta_\Kcal(2)$ is infinite}

In \cite{kriv2012} J.-L.~Krivine has shown that  $\Delta_\Kcal(2)$ does not contain 
any atoms (w.r.t.\ the order $\Delta_\Kcal(\leq_2)$), i.e.\ 
\[ \forall x{:}{\Delta_\Kcal(2)} 
   \bigl(x \neq 0  \to \exists y{:}{\Delta_\Kcal(2)}\, xy\neq 0 \wedge xy \neq x \bigr) \]
which by classical logic is equivalent to 
\[ \forall x{:}{\Delta_\Kcal(2)} 
\bigl( \forall y{:}{\Delta_\Kcal(2)}(xy \neq 0 \to xy \neq x \to \bot) \to x \neq 0  \to \bot\bigr) \]

For sake of completeness we recall Krivine's argument for which purpose we have to introduce
a bit of machinery. For $I \subseteq_\fin \N$ let $\bar{I} \in D$ with $\bar{I}(\vec{s}) = \top$ 
iff $s_i = \top_D$ for all $i \in I$. Notice that $\bar{\emptyset} = \top_D$ and
$\overline{\{n\}} = \bar{n}$. Obviously, we have 
\begin{enumerate}
\item[(1)] $\top_D \sqsubseteq u$ iff 
$u \in |\top,\bot \to \bot| \cap |\bot,\top \to \top|$
\item[(2)] $u \Vdash \bot,\bot \to \bot$ iff 
$u \in {\uparrow}\{\overline{I} \mid \emptyset \neq I \subseteq \{0,1\}\}$. 
\end{enumerate}
Let $t \in D$ with $t\top_D = \top_D$ and $t\overline{I} = \bar{0}$ for
nonempty subsets $I$ of $\{0,1\}$. Then $t$ realizes both
\[  |\top,\bot \to \bot| \cap |\bot,\top \to \top|, \top \to \bot 
    \qquad \mbox{and} \qquad (\bot,\bot \to \bot),\bot \to \bot \]
and thus $t$ realizes $\forall x{:}{\Delta_\Kcal(2)} 
\bigl( \forall y{:}{\Delta_\Kcal(2)}(xy \neq 0 \to xy \neq x \to \bot) \to x \neq 0  \to \bot\bigr)$
as can be seen by case analysis on $x \in \{0,1\}$.

Thus, in $\Kcal$ it holds that $\Delta_\Kcal(2)$ is infinite. But it is not clear 
\emph{a priori} whether $\Delta_\Kcal(2)$ is also \emph{Dedekind infinite}, i.e.\ 
whether the assertion
\[ \exists f{:}\Delta_\Kcal(2)^{N_\Kcal} (\forall n,m{:}N_\Kcal. f(n) \sim_{2_\Kcal} f(m) 
   \to n \sim_{N_\Kcal} m) \]
holds in $\Kcal$.\footnote{See e.g.\ \cite{jechAC} for the construction of a model of 
$\mathsf{ZF}$ in which there exists a Dedekind finite set which is not finite. This cannot be
achieved by forcing since forcing models all validate AC. One has to consider an appropriate
group $G$ of automorphism on an appropriate complete boolean algebra $B$ and take the 
$G$-invariant part of the $B$-valued model.} Actually, for quite some time we hoped that
in $\Kcal$ the object $\Delta_\Kcal(2)$ would not be Dedekind infinite since this would have
had the consequence that $\Kcal$ does not validate \emph{countable choice}.\footnote{As remarked
in \cite{jechAC} for any infinite set using countable choice one can prove quite straightforwardly
the existence of an injective function from $\N$ into this set.} The reason for this hope was that
presumably there does not exist a monomorphism $N_\Kcal \mono \Delta_\Kcal(2)$ in $\Kcal$. 

However, we will show that $\Kcal$ does indeed validate \emph{countable} and even 
\emph{dependent choice}. From this it follows that $\Kcal$ validates the assertion 
that there exists an injective function from $N_\Kcal$ to $\Delta_\Kcal(2)$ though presumably
this existential statement is not witnessed by a global element of $\Delta_\Kcal(2)^{N_\Kcal}$,
i.e.\ a proper monomorphism $N_\Kcal \mono \Delta_\Kcal(2)$ in $\Kcal$.

\section{$\Kcal$ validates countable and dependent choice}

Though Krivine's classical realizability gives rise to models of the classical 
set theory $\mathsf{ZF}$ (as described in \cite{kriv2001}) it generally does 
\emph{not} validate the full axiom of choice. Moreover, it is not known 
whether all classical realizability models for $\mathsf{ZF}$ validate the
principles of dependent or at least countable choice. Though, unfortunately, 
so far we do not know any counterexample J.-L.~Krivine strongly suspects that 
the answer to this question will be negative. In his opinion for realizing 
countable and dependent choice one has to extend his $\lambda$-calculus with 
control with new language constructs as described in \cite{kriv2003} where 
he adds a variant of LISP and Scheme's $\mathtt{quote}$ construct and shows 
how this may be used for realizing the above mentioned choice principles. But 
this method works only if the set $\Lambda$ of ``terms'' is countable which is,
obviously, not the case for the realizability structure arising from 
$D = \Sigma^{D^\omega}$ in $\mathbf{Coh}$ since $D$ has the size of the
continuum.

However, as known from work of C.~Spector dating back to the early 60s 
one may use \emph{bar recursion} for realizing classical choice principles.
This approach has been applied fruitfully in ``traditional'' proof theory
as described and discussed in U.~Kohlenbach's monograph \cite{ulrichbook}. 
However, Spector's original work and most of \cite{ulrichbook} are based on 
G\"odel's \emph{Dialectica} interpretation and not on realizability. Thus,
for our purposes the approach of U.~Berger and P.~Oliva in \cite{bo05}
(also discussed in \cite{ulrichbook}) is a better starting point since it 
is based on modified realizability which can be adapted more easily to the 
case of classical realizability.\footnote{This does not mean that methods
based on G\"odel's \emph{Dialectica} interpretation are not more appropriate
for the purposes of extracting programs and bounds from (classical) proofs
as emphasized in \cite{ulrichbook}.}\enlargethispage{\baselineskip}

In \cite{bo05} it is shown that when starting from a model $\Mo$ of 
higher type arithmetic validating an appropriate form of bar recursion
certain negative translations (where $\bot$ is replaced by arbitrary
$\Sigma^0_1$-formulas) of classical choice principles admit a modified 
realizability interpretation by objects of $\Mo$. For this purpose 
in \cite{bo05} they consider a `modified bar recursor' whose analogue 
in our setting we will introduce next after some preliminary remarks.

First of all for a coherence space $X$ we have to consider 
$X^* = \coprod\limits_{n\in\N} X^n$, the set of \emph{lists} of elements 
of $X$, which \emph{per se} is not a coherence space since it lacks a least 
element. However, we say that a map $f$ from $X^*$ to a coherence space $Y$ 
is \emph{stable} iff for all $n\in\N$ the restriction of $f$ to $X^n$ is 
stable. Moreover, if $Y$ is a coherence space then $X^*{\to}Y \cong
\prod_{n\in\omega} X^n{\to}Y$ is a coherence space since the $X^n{\to}Y$ are
coherence spaces. Alternatively, we may work in the slightly larger category
$\omega\mathbf{dI}_c$ of coherently complete countably based dI-domains
(see \cite{AmadioCurienBook}) and stable continuous functions between them.
We will have to consider stable functionals in the finite type hierarchy 
in $\mathbf{Coh}$ generated from $\Sigma$ and $D$ by $\to$, $(-)^\omega$ 
and $(-_1)^* \to (-_2)$. For every such type $X$ we have to specify its 
subset $\PL_X$ of \emph{proof-like} elements. 
Of course, for $D$ we put $\PL_D = P$ and for $\Sigma$ we put 
$\PL_\Sigma = \{\bot\}$. If $X$ and $Y$ are such types we put 
$\PL_{X{\to}Y} = \{ f : X \to Y \mid \forall x \in \PL_X. f(x) \in \PL_Y \}$,
$\PL_{X^\omega} = \PL_X^\omega$ and $\PL_{X^*} = \PL_X^*$, 
i.e.\ we extend $\PL$ \`a la logical relations.

\begin{defi}\label{BRdef}
Given  $Y : D^\omega \to \Sigma$ and $G : ((D \to \Sigma) \to \Sigma)^\omega$ 
in $\mathbf{Coh}$ let $\BR(Y,G)$ be the least stable function 
$\Psi : D^* \to \Sigma$ in $\mathbf{Coh}$ satisfying
\[ \Psi(s) = Y(s * \lambda n. G_{|s|}(\lambda x. \Psi(s * x))) \]
for all $s \in D^*$. 
\end{defi}
Obviously, the ensuing map 
$\BR : (D^\omega \to \Sigma) \to ((D \to \Sigma) \to \Sigma)^\omega\to \Sigma$ 
is stable and proof-like.

Notice that all the types built from $D$ and $\Sigma$ by $\to$, 
$(-)^\omega$ and $(-_1)^* \to (-_2)$ appear as retracts of $D$ via 
proof-like maps.
They form a typed pca realizability over which gives rise to a category
equivalent to $\RT(D,P)$ as described on a more general level 
in \cite{LietzStreicher}. 
This allows us to assume that realizers of particular propositions have 
particular types which often allows us to reason 
in a more intuitive way.

We often will have to refer to $\Omega_\Kcal$ considered as an object of $\E$.
This object has underlying set $\{ A \in \Pow(D) \mid A^{\Perp\Perp} = A \}$ for
which equality is given by logical equivalence\footnote{which is the same in $\E$ 
and $\Kcal$ for propositions of this particular form}. Moreover, for objects $X$
in $\Kcal$ the exponential $\Pow_\Kcal(X) = \Omega_\Kcal^X$ is the same when taken
in $\E$ and $\Kcal$, respectively. Moreover, for $X$ in $\E$ the map
$\Omega_\Kcal^{\eta_X} : \Omega_\Kcal^{i_*i^*X} \to \Omega_\Kcal^X$ (where 
$\eta_X : X \to i_*i^*X$ is the unit of $i^* \dashv i_*$ at $X$) is an isomorphism 
in $\E$. Accordingly, we will often write $\Pow_\Kcal(X)$ for $\Omega_\Kcal^X$ 
in $\E$.

\subsection{Countable Choice holds in $\Kcal$}

Spector already observed that the negative translation of countable choice 
can be proved in any intuitionistic theory validating  countable choice and
the principle of \emph{Double Negation Shift} (DNS) for formulas in the 
negative fragment. Like all relative realizability toposes $\E = \RT(D,P)$ 
validates countable and dependent choice. Thus, due to Spector's observation
it suffices to show that $\E$ also validates an appropriate form of DNS.

\begin{lem}\label{DNSE}
The topos $\E$ validates the principle
\[ \mathrm{(DNS)} \qquad  \forall B{:}\Pow_\Kcal(N). 
  \forall n. {\sim\sim}B(n) \to \;{\sim\sim}\forall n.B(n) \]
where ${\sim}A$ stands for $A \to U$.
\end{lem}
\proof
Suppose $B \in \Pow(D)^\omega$ with $B(n)^{\Perp\Perp} = B(n)$ for all $n$,
$G$ of type $((D \to \Sigma) \to \Sigma)^\omega$ realize 
$\forall n.{\sim\sim}B(n)$ and $Y$ of type $D^\omega \to \Sigma$ 
realize ${\sim}\forall n.B(n)$. Let $\Psi = \BR(Y,G)$. Using a variant of 
bar induction as described in \cite{bo05} we will show now that 
$\Psi(\langle\rangle) = \top$ and thus realizes $U$.

We write $S(x,n)$ for $x \in B(n)$ and $P(s)$ for $\Psi(s) = \top$. 
We employ the abbreviations $s \in S \equiv \forall k < |s|\; s_k \in B(k)$ 
and $\alpha \in S \equiv \forall k\; \alpha_k \in B(k)$. 
By bar induction relativized to $S$ (see \cite{bo05} for details) 
for showing $P(\langle\rangle)$ it suffices to show that
\begin{enumerate}
\item[(1)] $\forall \alpha \in S \exists n \; P(\bar{\alpha}(n))$
\item[(2)] $\forall s \in S \bigl(\forall x (S(x,|s|) \to P(s*x))\bigr) \to P(s)$.
\end{enumerate}
\emph{ad} (1) : Suppose $\alpha \in S$, i.e.\ $\alpha(n) \in B(n)$ for all $n$.
Then by assumption on $Y$ we have $Y(\alpha) = \top$. Since $Y$ is continuous
there exists an $n$ with $Y(\alpha) = Y(\bar{\alpha}(n)*\beta)$ for all $\beta$.
Thus, we have $\Psi(\bar{\alpha}(n)) = \top$, i.e.\ $P(\bar{\alpha}(n))$ as 
desired.

\noindent
\emph{ad} (2) : 
Suppose $s \in S$ with $\forall x (S(x,|s|) \to P(s{*}x))$, i.e.\
$\forall x (x \in B(|s|) \to \Psi(s{*}x) = \top )$. Thus 
$\lambda x. \Psi(s*x)$ realizes ${\sim}B(|s|)$. Accordingly, by assumption 
on $G$ it follows that $G_{|s|}(\lambda x. \psi(s*x))$ realizes $U$ and thus 
also $B(n)$ (since $B(n)$ contains $U$ as a subset like all propositions 
in $\Kcal$). Thus $s * \lambda n. G_{|s|}(\lambda x. \psi(s*x))$ realizes
$\forall n. B(n)$ and, accordingly, by assumption on $Y$ it follows that
$\Psi(s) = Y(s * \lambda n. G_{|s|}(\lambda x. \psi(s*x)))$ realizes $U$,
i.e.\ $P(s)$ as desired.

Thus, since $\lambda G. \lambda Y. \BR(Y,G)(\langle\rangle)$ is proof-like 
it realizes the proposition 
\[  \forall B{:}\Pow_\Kcal(N). 
    \forall n. {\sim\sim}B(n) \to \;{\sim\sim}\forall n.B(n) \]
which, therefore, holds in $\E$ as claimed.
\qed

Notice that the form of bar induction used in the proof of Th.~\ref{DNSE}
is valid only because $D^\omega$ consists of {\bf all} sequences in $D$ 
(and not just the computable ones).

Now we are ready to show that countable choice holds in $\Kcal$.

\begin{thm}\label{CCK}
For every object $X$ in $\Kcal$ the proposition
\[ \forall R{:}\Pow(N{\times}X).\, \forall n{:}N.\exists x{:}X. R(n,x) \to
                                   \exists f{:}X^N.\forall n{:}N.R(n,f(n)) \]
hold in $\Kcal$.
\end{thm}
\proof
Since $\Kcal$ is equivalent to the subtopos $\E_U$ of $\E$ consisting of sheaves 
for $j_U = \;\sim\circ\sim$ the problem reduces to showing that
\[ \forall R{:}\Pow_\Kcal(N{\times}X). \, 
          \forall n{:}N.{\sim\sim}\exists x{:}X. R(n,x) \to
                {\sim\sim}\exists f{:}X^N.\forall n{:}N.R(n,f(n)) \]
holds in $\E$. By Lemma~\ref{DNSE} the implication
\[ \forall n{:}N.{\sim\sim}\exists x{:}X. R(n,x) \to 
    {\sim\sim}\forall n{:}N.\exists x{:}X. R(n,x) \]
holds in $\E$ and thus it suffices to show that
\[ \forall R{:}\Pow_\Kcal(N{\times}X). \, 
          {\sim\sim}\forall n{:}N.\exists x{:}X. R(n,x) \to
                {\sim\sim}\exists f{:}X^N.\forall n{:}N.R(n,f(n)) \]
holds in $\E$. This, however, holds since $\E$ validates countable choice
and $\sim\sim$ commutes with implication.
\qed

Thus, we have shown that $\Kcal$ validates countable choice since $X^N$ is
isomorphic to $X^{N_\Kcal}$ in $\E$.

\smallskip
Notice that for classical realizability models arising from
countable term models one cannot apply the method we have used here
because bar induction does not seem to be applicable since not every 
external sequence of terms can be represented by a term. Thus, for
countable term models Krivine in \cite{kriv2003} introduced a
$\mathtt{quote}$-like construct for the purpose of realizing
countable choice. Apparently, these two different methods are 
applicable under \emph{mutually exclusive circumstances}. Whether
countable choice holds in all realizability models is unknown up to
now but one strongly suspects that the answer is negative!

\subsection{Dependent Choice in $\Kcal$}

A topos with natural numbers object $N$ validates the principle 
DC of \emph{Dependent Choice} iff 
\begin{tabbing}
\qquad \= $\forall R{:}\Pow(N{\times}X{\times}X).$\=
           $\forall n{:}N.\forall x{:}X.\exists y{:}X. R(n,x,y) \to$\\
\>\> \quad $\forall a{:}X.\exists f{:}X^N.\, f(0) = a \wedge 
      \forall n{:}N. R(n,f(n),f(n{+}1))$
\end{tabbing}
holds for every object $X$ of the topos. It is well know that $\E$ and actually 
every relative realizability topos validates DC. Unfortunately,
the validity of Double Negation Shift in $\E$ is not sufficient for reducing
validity of DC in $\Kcal$ to its validity in $\E$. For this reason 
in Theorem~4 of \cite{bo05} it is shown how to use modified bar recursion
for realizing appropriate negative translations of DC. With some effort 
their proof can be adapted to $\Kcal$. We leave the tedious details to the
inclined reader. Notice, however, that Theorem~\ref{CCK} suffices already
for showing that the infinite object $\Delta_\Kcal(2)$ is also Dedekind 
infinite, i.e.\ that $\Kcal$ validates the proposition
$\exists f{:}\Delta_\Kcal(2)^N (\forall n,m{:}N. f(n) \sim_{2_\Kcal} f(m)$.
However, this valid existential statement need not be witnessed by a global
element of $\Delta_\Kcal(2)^N$.

\section{Is $\Kcal$ 2-valued?}\label{isK2val}

A proposition $A \in \Omega_\Kcal$ is \emph{valid in $\Kcal$} iff 
$A \cap P \neq \emptyset$. The topos $\Kcal$ is $2$-valued iff for every
$A \in \Omega_\Kcal$ either $A$ or $\neg A$ has nonempty intersection with $P$.

Notice that for $t \in D$ we have $t \in P$ iff $\trace(t) \cap P^\omega =
\emptyset$. 
Thus, if $A$ holds in $\Kcal$ then $A^\Perp \cap P^\omega = \emptyset$. If the 
reverse implication held as well then $\Kcal$ would be $2$-valued which can be 
seen as follows. 
Suppose $A$ does not hold in $\Kcal$. Then, due to our assumption, there exists 
$\vec{s} \in A^\Perp \cap P^\omega$ and thus $\lambda t.\vec{s} \in D^\omega.\,
t(\vec{s})$ is an element of $P \cap \neg A$. 

But if $A$ is the biorthogonal closure of a countable subset of $D$ 
we actually can reverse the implication.

\begin{lem}
If $A = \{ t_n \mid n\in\omega\}^{\Perp\Perp}$ 
with $A^\Perp \cap P^\omega = \emptyset$ then $A \cap P \neq \emptyset$.
\end{lem}
\proof
W.l.o.g.\footnote{This can be achieved easily since 
$\sqcap : \Sigma \times \Sigma \to \Sigma$ is stable.} we assume that 
$t_{n+1}^{-1}(\top) \subseteq t_n^{-1}(\top)$ for all $n \in \omega$.
We consider the countably branching tree $T = \bigcup\limits_{n\in\omega}
\{n\} \times \trace(t_n)$ where the ancestor of $\pair{n{+}1}{\vec{s}}$
is the unique element $\pair{n}{\vec{r}}$ with $\vec{r} \sqsubseteq \vec{s}$.
Observe that for every $\vec{s} \in A^\Perp$ and $n\in\omega$ there is a
unique $\vec{s}^{(n)} \in \trace(t_n)$ with $\vec{s}^{(n)} \sqsubseteq \vec{s}$.
Thus, the minimal elements of $A^\Perp$ are precisely the suprema of the
infinite paths in $T$, i.e.\ for every $\vec{s} \in \min(A^\Perp)$ we have
$\vec{s} = \bigsqcup\limits_{n\in\omega} \vec{s}^{(n)}$. Thus, due to our
assumption $A^\Perp \cap P^\omega = \emptyset$ every infinite path through $T$ 
eventually leads out of $P^\omega$. Let $t$ be the element of $D$ whose 
trace consists of those finite elements $\vec{s}$ of $D^\omega$ with 
$\vec{s}^{(n)} \not\in P^\omega$ but $\vec{s}^{(k)} \in P^\omega$ for all $k < n$. 
Obviously, we have $t \in P$ and $\min(A^\Perp) \subseteq t^{-1}(\top)$. 
Thus $t \in A \cap P$ as desired.
\qed

In order to generalize this lemma to arbitrary propositions in $\Kcal$ one 
could try to work with a well ordering of a biorthogonally closed subset $A$ 
of $D$ but then beyond stage $\omega$ the labels of the tree $T$ are not 
finite anymore. 

Another line of attack would be as follows. Suppose $A = A^{\Perp\Perp}$ such
that $\trace(t) \cap P^\omega \neq\emptyset$ for all $t \in A$. Notice that
the (upward closures) of the sets $\trace(t) \cap P^\omega$ with $t \in A$ form 
a filter w.r.t.\ the Smyth ordering. But, alas, we do not know how to prove
that the intersection of the elements of this filter has to be non-empty.

On the other hand we do not know any particular biorthogonally closed subset 
of $D$ which does not already arise as the biorthogonal of a countable subset. 
In particular, we may replace any proposition $A$ with the biorthogonal
closure of the intersection of $A$ with the computable elements of $A$.
Maybe this does not make any difference for propositions $A$ arising from
the interpretation of a closed formula in the language of set theory.
 
\section{Summary}

We have shown that a new boolean non-Grothendieck topos $\Kcal$
arises from a canonical model of $\lambda$-calculus with control in
the category $\mathbf{Coh}$ of coherence spaces and stable functions.
We have shown that $\Kcal$ validates all true sentences of first order 
arithmetic and also countable (and dependent) choice. 

We have also observed that the model constructions collapses to the 
ground model $\Set$ when starting from the canonical model of 
$\lambda$-calulus with control in Scott domains where as usual the 
culprit is parallel-or. 

There are still quite a few open questions about the topos $\Kcal$ arising
from the stable model of $\lambda$-calulus with control. One would like to
see a concrete example of a set-theoretic statement holding in $\Set$ but not
in $\Kcal$. We suspect that AC, the full axiom of choice, is such an example 
but have not been able yet to verify this. Moreover, one would like to know 
whether every closed formula in the language of set theory is decided 
by $\Kcal$.

All our arguments apply also to the case when $D \cong \Sigma^{D^\omega}$ is
solved in the category of observably sequential algorithms because it does
not host parallel-or and all (even non-effective) elements of $D^\omega$ are 
represented by elements of $D$.

\section*{Acknowledgement}

I want to thank J.-L.~Krivine for patiently explaining to me the intuitions
underlying his work on classical realizability.


\end{document}